\documentclass[11pt, twoside]{article}
\usepackage{latexsym}
\usepackage{amsmath}
\usepackage{amssymb}
\usepackage[all]{xy}
\usepackage{amsfonts}
\usepackage{verbatim}
\usepackage{amsthm}
\usepackage{mathrsfs}
\usepackage{epsfig}
\usepackage{xy}
\usepackage{array}
\usepackage{stmaryrd}
\usepackage{graphicx,color}
\usepackage{xcolor}
\usepackage[colorlinks=true,linkcolor=blue,citecolor=blue]{hyperref}
\usepackage{tikz}
\usetikzlibrary{arrows,calc}
\usepackage{etex}
\usepackage{mathdots}
\usepackage{float}
\usepackage{graphics}
\usepackage{pdflscape}
\usepackage{extarrows}
\usepackage{anysize,hyperref}
\input xypic
\xyoption{all}
\usepackage{bm}
\usepackage[perpage,symbol]{footmisc}
\usepackage{setspace}
\usepackage{enumitem}

\def\A{\mathcal{A}}
\def\T{\mathcal{T}}

\def\C{\mathscr{C}}
\def\E{\mathbb{E}}

\def\s{\mathfrak{s}}
\def\id{\mathrm{id}}
\def\op{^\mathrm{op}}
\def\Ab{\mathsf{Ab}}
\def\del{\delta}
\def\rad{\mbox{rad}}
\def\dr{\ar@{->}[r]}

\def\X{\mathscr{X}}

\def\add{\mbox{add}}

\def\Hom{\mbox{Hom}}

\def\gh{\mathsf{Gh}\hspace{.01in}}
\def\cogh{\mathsf{CoGh}\hspace{.01in}}

\newcommand{\CC}{{\bf{C}}^{n+2}_{\C}}
\newcommand{\mr}{\hbox{\boldmath$\cdot$}}
\newcommand{\ov}{\overset}
\newcommand{\lra}{\longrightarrow}
\newcommand{\co}{\colon}
\newcommand{\uas}{^{\ast}}            
\newcommand{\sas}{_{\ast}}
\newcommand{\Xd}{\langle X^{\mr},\del\rangle}  

\newcommand{\ush}{^\sharp}           
\newcommand{\ssh}{_\sharp}
\SelectTips{cm}{10}

\def\G{\mathbb{G}}
\def\H{\mathbb{H}}
\begin{document}
\title{\Large{\bf Pre-$\bm{(n+2)}$-angulated categories$^\bigstar$\footnotetext{\hspace{-1em}$^\bigstar$Panyue Zhou was supported by the National Natural Science Foundation of China (Grant No. 11901190).}}}
\medskip
\author{Jing He, Panyue Zhou and Xingjia Zhou}

\date{}

\maketitle
\def\blue{\color{blue}}
\def\red{\color{red}}

\newtheorem{theorem}{Theorem}[section]
\newtheorem{lemma}[theorem]{Lemma}
\newtheorem{corollary}[theorem]{Corollary}
\newtheorem{proposition}[theorem]{Proposition}
\newtheorem{conjecture}{Conjecture}
\theoremstyle{definition}
\newtheorem{definition}[theorem]{Definition}
\newtheorem{question}[theorem]{Question}
\newtheorem{remark}[theorem]{Remark}
\newtheorem{remark*}[]{Remark}
\newtheorem{example}[theorem]{Example}
\newtheorem{example*}[]{Example}
\newtheorem{condition}[theorem]{Condition}
\newtheorem{condition*}[]{Condition}
\newtheorem{construction}[theorem]{Construction}
\newtheorem{construction*}[]{Construction}

\newtheorem{assumption}[theorem]{Assumption}
\newtheorem{assumption*}[]{Assumption}

\baselineskip=17pt
\parindent=0.5cm
\vspace{-6mm}

\begin{abstract}
\baselineskip=16pt
In this article, we introduce the notion of pre-$(n+2)$-angulated categories
as higher dimensional analogues of pre-triangulated categories defined by Beligiannis-Reiten.
We first show that the idempotent completion of a pre-$(n+2)$-angulated category admits a unique structure
of pre-$(n+2)$-angulated category.
Let $(\C,\E,\s)$ be an $n$-exangulated category and
$\X$ be a strongly functorially finite subcategory of $\C$.
We then show that the quotient category $\C/\X$ is a pre-$(n+2)$-angulated category.
These results allow to construct several examples of pre-$(n+2)$-angulated categories.
Moreover, we also give a necessary and sufficient
condition for the quotient $\C/\X$ to be an $(n+2)$-angulated
category.\\[0.2cm]
\end{abstract}

\pagestyle{myheadings}
\markboth{\rightline {\scriptsize Jing He,  Panyue Zhou and Xingjia Zhou}}
         {\leftline{\scriptsize Pre-$(n+2)$-angulated categories}}

\section{Introduction}
Triangulated categories were introduced in the mid 1960's by Verdier \cite{V}
Having their origins in algebraic geometry and algebraic topology, triangulated categories have by now become indispensable in many different areas of mathematics.
Geiss, Keller and Oppermann \cite{GKO} introduced
a new type of categories, called $(n+2)$-angulated categories, which generalize triangulated categories: the classical triangulated categories are the special case $n=1$.
Jasso \cite{J} introduced $n$-abelian and $n$-exact categories, these are analogs of abelian and
exact categories from the point of view of higher homological algebra.
The case $n=1$ corresponds to the usual concepts of abelian and exact categories. Later, Lin \cite{L2} introduced right $(n+2)$-angulated categories, the right triangulated categories in the sense of Assem, Beligiannis and Marmaridis \cite{ABM} are the special case $n=1$.

Extriangulated categories were recently introduced by Nakaoka and Palu \cite{NP} by extracting those properties of ${\rm Ext}^1$ on exact categories and on triangulated categories that seem relevant from the point of view of cotorsion pairs. In particular, triangulated categories and exact categories are extriangulated categories. There are many examples of extriangulated categories which are neither triangulated categories nor exact categories, see
see \cite{NP,ZZ,HZZ,ZhZ,NP1}.
Recently, Herschend, Liu and Nakaoka \cite{HLN1} introduced the notion of $n$-exangulated categories. It should be noted that the case $n =1$ corresponds to extriangulated categories. As typical examples we have that $n$-exact and $(n+2)$-angulated categories are $n$-exangulated, see \cite[Proposition 4.5 and Proposition 4.34]{HLN1}. However, there are some other examples of $n$-exangulated categories which are neither $(n+2)$-angulated nor $n$-exact, see \cite{HLN1,HLN2,LZ,HZZ,HZZ1,HHZ,KMS}.

Beligiannis and Reiten \cite{BR} introduced the notion of pre-triangulated categories
which provides a common generalization of abelian,
triangulated and stable categories.
Based on Beligiannis-Reiten's idea, we introduce the notion of
pre-$(n+2)$-angulated categories, see Definition \ref{pre}.
Our first main result is the following.

\begin{theorem}\label{main111}{\rm (see Theorem \ref{main0} for more details)}
Let $(\C,\Sigma^n,\Omega^n,\Theta,\Phi,\varepsilon,\eta)$ be a pre-$(n+2)$-angulated category.
Then its idempotent completion $(\widetilde{\C},\widetilde{\Sigma^n},\widetilde{\Omega^n},\widetilde{\Theta},\widetilde{\Phi},\widetilde{\varepsilon},\widetilde{\eta})$  admits a unique structure of pre-$(n+2)$-angulated category.
\end{theorem}

Our second result shows
how to equip quotient categories with a pre-$(n+2)$-angulated structure as follows.

\begin{theorem}\label{main3}
{\rm (see Theorem \ref{main1} for more details)}
Let $(\C,\E,\s)$ be an $n$-exangulated category and
$\X$ be a strongly functorially finite subcategory of $\C$.
Then the quotient category $\C/\X$ is a pre-$(n+2)$-angulated category.
\end{theorem}

As an application of Theorem \ref{main3}, we show
how to equip quotient categories with an $(n+2)$-angulated structure.
This recovers and gives a short proof to a result of Zhou \cite[Theorem 5.7]{Z}.

\begin{theorem}\label{main4}{\rm (see Theorem \ref{main6} for more details)}
Let $(\C,\Sigma^n,\Theta)$ be an $(n+2)$-angulated category with
 a Serre functor $\mathbb{S}$, and $\X$ be a strongly functorially finite subcategory of $\C$.
 Then the pre-$(n+2)$-angulated category $\C/\X$ is $(n+2)$-angulated if and only $\mathbb{S}\X=\Sigma^n\X$.
\end{theorem}

This article is organized as follows. In Section 2, we review some elementary definitions and facts on $n$-exangulated categories. In Section 3, we define the notion of pre-$(n+2)$-angulated categories
and give some examples of pre-$(n+2)$-angulated categories. In Section 4, we give the proof of Theorem \ref{main111}.
In Section 5, we give the proof of Theorem \ref{main3}. In Section 6, we give the proof of Theorem \ref{main4}.

\section{Preliminaries}
In this section, we briefly review basic concepts and results on $n$-exangulated categories.

 Let $\C$ be an additive category and
 $\mathbb{E}\colon \C^{\rm op}\times \C \rightarrow \Ab~~\mbox{($\Ab$ is the category of abelian groups)}$
 an additive bifunctor.
 {For any pair of objects $A,C\in\C$, an element $\del\in\E(C,A)$ is called an {\it $\E$-extension} or simply an {\it extension}. We also write such $\del$ as ${}_A\del_C$ when we indicate $A$ and $C$. The zero element ${}_A0_C=0\in\E(C,A)$ is called the {\it split $\E$-extension}. For any pair of $\E$-extensions ${}_A\del_C$ and ${}_{A'}\del{'}_{C'}$, let $\delta\oplus \delta'\in\mathbb{E}(C\oplus C', A\oplus A')$ be the
element corresponding to $(\delta,0,0,{\delta}{'})$ through the natural isomorphism $\mathbb{E}(C\oplus C', A\oplus A')\simeq\mathbb{E}(C, A)\oplus\mathbb{E}(C, A')
\oplus\mathbb{E}(C', A)\oplus\mathbb{E}(C', A')$.

For any $a\in\C(A,A')$ and $c\in\C(C',C)$,  $\E(C,a)(\del)\in\E(C,A')\ \ \text{and}\ \ \E(c,A)(\del)\in\E(C',A)$ are simply denoted by $a_{\ast}\del$ and $c^{\ast}\del$, respectively.

Let ${}_A\del_C$ and ${}_{A'}\del{'}_{C'}$ be any pair of $\E$-extensions. A {\it morphism} $(a,c)\colon\del\to{\delta}{'}$ of extensions is a pair of morphisms $a\in\C(A,A')$ and $c\in\C(C,C')$ in $\C$, satisfying the equality
$a_{\ast}\del=c^{\ast}{\delta}{'}$.}

Let $\C$ be an additive category and $n$ be any positive integer.
\begin{definition}\cite[Definition 2.7]{HLN1}
Let $\bf{C}_{\C}$ be the category of complexes in $\C$. As its full subcategory, define $\CC$ to be the category of complexes in $\C$ whose components are zero in the degrees outside of $\{0,1,\ldots,n+1\}$. Namely, an object in $\CC$ is a complex $X_{\bullet}=\{X_i,d^X_i\}$ of the form
\[ X_0\xrightarrow{d^X_0}X_1\xrightarrow{d^X_1}\cdots\xrightarrow{d^X_{n-1}}X_n\xrightarrow{d^X_n}X_{n+1}. \]
We write a morphism $f_{\bullet}\co X_{\bullet}\to Y_{\bullet}$ simply $f_{\bullet}=(f_0,f_1,\ldots,f_{n+1})$, only indicating the terms of degrees $0,\ldots,n+1$.
\end{definition}

\begin{definition}\cite[Definition 2.11]{HLN1}
By Yoneda lemma, any extension $\del\in\E(C,A)$ induces natural transformations
\[ \del\ssh\colon\C(-,C)\Rightarrow\E(-,A)\ \ \text{and}\ \ \del\ush\colon\C(A,-)\Rightarrow\E(C,-). \]
For any $X\in\C$, these $(\del\ssh)_X$ and $\del\ush_X$ are given as follows.
\begin{enumerate}
\item[\rm(1)] $(\del\ssh)_X\colon\C(X,C)\to\E(X,A)\ :\ f\mapsto f\uas\del$.
\item[\rm (2)] $\del\ush_X\colon\C(A,X)\to\E(C,X)\ :\ g\mapsto g\sas\delta$.
\end{enumerate}
We simply denote $(\del\ssh)_X(f)$ and $\del\ush_X(g)$ by $\del\ssh(f)$ and $\del\ush(g)$, respectively.
\end{definition}

\begin{definition}\cite[Definition 2.9]{HLN1}
 Let $\C,\E,n$ be as before. Define a category $\AE:=\AE^{n+2}_{(\C,\E)}$ as follows.
\begin{enumerate}
\item[\rm(1)] A object in $\AE^{n+2}_{(\C,\E)}$ is a pair $\Xd$ of $X_{\bullet}\in\CC$
and $\del\in\E(X_{n+1},X_0)$ satisfying
$$(d_0^X)_{\ast}\del=0~~\textrm{and}~~(d^X_n)^{\ast}\del=0.$$
We call such a pair an $\E$-attached
complex of length $n+2$. We also denote it by
$$X_0\xrightarrow{d_0^X}X_1\xrightarrow{d_1^X}\cdots\xrightarrow{d_{n-2}^X}X_{n-1}
\xrightarrow{d_{n-1}^X}X_n\xrightarrow{d_n^X}X_{n+1}\overset{\delta}{\dashrightarrow}.$$
\item[\rm (2)]  For such pairs $\Xd$ and $\langle Y_{\bullet},\rho\rangle$, a morphism $f_{\bullet}\colon\Xd\to\langle Y_{\bullet},\rho\rangle$ is
defined to be a morphism $f_{\bullet}\in\CC(X_{\bullet},Y_{\bullet})$ satisfying $(f_0)_{\ast}\del=(f_{n+1})^{\ast}\rho$.

We use the same composition and the identities as in $\CC$.

\end{enumerate}
\end{definition}

\begin{definition}\cite[Definition 2.13]{HLN1}\label{def1}
 An {\it $n$-exangle} is a pair $\Xd$ of $X_{\bullet}\in\CC$
and $\del\in\E(X_{n+1},X_0)$ which satisfies the following conditions.
\begin{enumerate}
\item[\rm (1)] The following sequence of functors $\C\op\to\Ab$ is exact.
$$
\C(-,X_0)\xrightarrow{\C(-,\ d^X_0)}\cdots\xrightarrow{\C(-,\ d^X_n)}\C(-,X_{n+1})\xrightarrow{~\del\ssh~}\E(-,X_0)
$$
\item[\rm (2)] The following sequence of functors $\C\to\Ab$ is exact.
$$
\C(X_{n+1},-)\xrightarrow{\C(d^X_n,\ -)}\cdots\xrightarrow{\C(d^X_0,\ -)}\C(X_0,-)\xrightarrow{~\del\ush~}\E(X_{n+1},-)
$$
\end{enumerate}
In particular any $n$-exangle is an object in $\AE$.
A {\it morphism of $n$-exangles} simply means a morphism in $\AE$. Thus $n$-exangles form a full subcategory of $\AE$.
\end{definition}

\begin{definition}\cite[Definition 2.22]{HLN1}
Let $\s$ be a correspondence which associates a homotopic equivalence class $\s(\del)=[{}_A{X_{\bullet}}_C]$ to each extension $\del={}_A\del_C$. Such $\s$ is called a {\it realization} of $\E$ if it satisfies the following condition for any $\s(\del)=[X_{\bullet}]$ and any $\s(\rho)=[Y_{\bullet}]$.
\begin{itemize}
\item[{\rm (R0)}] For any morphism of extensions $(a,c)\co\del\to\rho$, there exists a morphism $f_{\bullet}\in\CC(X_{\bullet},Y_{\bullet})$ of the form $f_{\bullet}=(a,f_1,\ldots,f_n,c)$. Such $f_{\bullet}$ is called a {\it lift} of $(a,c)$.
\end{itemize}
In such a case, we simply say that \lq\lq$X_{\bullet}$ realizes $\del$" whenever they satisfy $\s(\del)=[X_{\bullet}]$.

Moreover, a realization $\s$ of $\E$ is said to be {\it exact} if it satisfies the following conditions.
\begin{itemize}
\item[{\rm (R1)}] For any $\s(\del)=[X_{\bullet}]$, the pair $\Xd$ is an $n$-exangle.
\item[{\rm (R2)}] For any $A\in\C$, the zero element ${}_A0_0=0\in\E(0,A)$ satisfies
\[ \s({}_A0_0)=[A\ov{\id_A}{\lra}A\to0\to\cdots\to0\to0]. \]
Dually, $\s({}_00_A)=[0\to0\to\cdots\to0\to A\ov{\id_A}{\lra}A]$ holds for any $A\in\C$.
\end{itemize}
Note that the above condition {\rm (R1)} does not depend on representatives of the class $[X_{\bullet}]$.
\end{definition}

\begin{definition}\cite[Definition 2.23]{HLN1}
Let $\s$ be an exact realization of $\E$.
\begin{enumerate}
\item[\rm (1)] An $n$-exangle $\Xd$ is called an $\s$-{\it distinguished} $n$-exangle if it satisfies $\s(\del)=[X_{\bullet}]$. We often simply say {\it distinguished $n$-exangle} when $\s$ is clear from the context.
\item[\rm (2)]  An object $X_{\bullet}\in\CC$ is called an {\it $\s$-conflation} or simply a {\it conflation} if it realizes some extension $\del\in\E(X_{n+1},X_0)$.
\item[\rm (3)]  A morphism $f$ in $\C$ is called an {\it $\s$-inflation} or simply an {\it inflation} if it admits some conflation $X_{\bullet}\in\CC$ satisfying $d_0^X=f$.
\item[\rm (4)]  A morphism $g$ in $\C$ is called an {\it $\s$-deflation} or simply a {\it deflation} if it admits some conflation $X_{\bullet}\in\CC$ satisfying $d_n^X=g$.
\end{enumerate}
\end{definition}

\begin{definition}\cite[Definition 2.27]{HLN1}
For a morphism $f_{\bullet}\in\CC(X_{\bullet},Y_{\bullet})$ satisfying $f_0=\id_A$ for some $A=X_0=Y_0$, its {\it mapping cone} $M_{_{\bullet}}^f\in\CC$ is defined to be the complex
\[ X_1\xrightarrow{d^{M_f}_0}X_2\oplus Y_1\xrightarrow{d^{M_f}_1}X_3\oplus Y_2\xrightarrow{d^{M_f}_2}\cdots\xrightarrow{d^{M_f}_{n-1}}X_{n+1}\oplus Y_n\xrightarrow{d^{M_f}_n}Y_{n+1} \]
where $d^{M_f}_0=\begin{bmatrix}-d^X_1\\ f_1\end{bmatrix},$
$d^{M_f}_i=\begin{bmatrix}-d^X_{i+1}&0\\ f_{i+1}&d^Y_i\end{bmatrix}\ (1\le i\le n-1),$
$d^{M_f}_n=\begin{bmatrix}f_{n+1}&d^Y_n\end{bmatrix}$.

{\it The mapping cocone} is defined dually, for morphisms $h_{\bullet}$ in $\CC$ satisfying $h_{n+1}=\id$.
\end{definition}

\begin{definition}\cite[Definition 2.32]{HLN1}
An {\it $n$-exangulated category} is a triplet $(\C,\E,\s)$ of additive category $\C$, additive bifunctor $\E\co\C\op\times\C\to\Ab$, and its exact realization $\s$, satisfying the following conditions.
\begin{itemize}[leftmargin=3.3em]
\item[{\rm (EA1)}] Let $A\xrightarrow{f}B\xrightarrow{g}C$ be any sequence of morphisms in $\C$. If both $f$ and $g$ are inflations, then so is $g\circ f$. Dually, if $f$ and $g$ are deflations, then so is $g\circ f$.
\end{itemize}
\vspace{-0.7cm}

\begin{itemize}[leftmargin=3.3em]
\item[{\rm (EA2)}] For $\rho\in\E(D,A)$ and $c\in\C(C,D)$, let ${}_A\langle X_{\bullet},c\uas\rho\rangle_C$ and ${}_A\langle Y_{\bullet},\rho\rangle_D$ be distinguished $n$-exangles. Then $(1_A,c)$ has a {\it good lift} $f_{\bullet}$, in the sense that its mapping cone gives a distinguished $n$-exangle $\langle M^f_{_{\bullet}},(d^X_0)\sas\rho\rangle$.
\end{itemize}
\vspace{-0.7cm}

\begin{itemize}[leftmargin=4.2em]
 \item[{\rm (EA2$\op$)}] Dual of {\rm (EA2)}.
\end{itemize}
 \vspace{-0.3cm}

 Note that in the case $n=1$, a triplet $(\C,\E,\s)$ is a  $1$-exangulated category if and only if it is an extriangulated category, see \cite[Proposition 4.3]{HLN1}.
\end{definition}

\begin{example}
From \cite[Proposition 4.34]{HLN1} and \cite[Proposition 4.5]{HLN1},  we know that $n$-exact categories and $(n+2)$-angulated categories are $n$-exangulated categories.
There are some other examples of $n$-exangulated categories
 which are neither $n$-exact nor $(n+2)$-angulated, see \cite{HLN1,HLN2,LZ,HZZ1}.
\end{example}

\section{Pre-$(n+2)$-angulated categories}
In this section we introduce the concept of a pre-$(n+2)$-angulated category and give  many examples of pre-$(n+2)$-angulated categories.

If an additive category is left and right $(n+2)$-angulated,
then usually the left and right  $(n+2)$-angulated structures are compatible in a nice way.
Now we formalize this situation in the following definition.

\begin{definition}\label{pre}
Let $\C$ be an additive category. A pre-$(n+2)$-angulation of $\C$
consists of the following data:
\begin{itemize}
\item[(i)] An adjoint pair $(\Sigma^n,\Omega^n)$ of additive endofunctors $\Sigma^n,\Omega^n\colon
\C\to \C$. Let $\varepsilon\colon \Sigma^n\Omega^n\to 1_{\C}$
be the counit and let $\eta\colon 1_{\C}\to \Omega^n\Sigma^n$
be the unit of the adjoint pair.

\item[(ii)] A collection of diagrams $\Phi$ in $\C$ of the form
$$\Omega^nA_{n+1}\xrightarrow{~}A_0\xrightarrow{~}A_1\xrightarrow{~}A_2\xrightarrow{~}\cdots\xrightarrow{~}A_n\xrightarrow{~}A_{n+1}$$
such that the triple $(\C,\Omega^n,\Phi)$  is a left $(n+2)$-angulated category.

\item[(iii)] A collection of diagrams $\Theta$ in $\C$ of the form
$$A_0\xrightarrow{~}A_1\xrightarrow{~}A_2\xrightarrow{~}\cdots\xrightarrow{~}A_n\xrightarrow{~}A_{n+1}
\xrightarrow{~}\Sigma^nA_0$$
such that the triple $(\C,\Sigma^n,\Theta)$  is a right $(n+2)$-angulated category.

\item[(iv)] Each solid commutative diagram in $\C$:
$$\xymatrix{
A_0 \ar[r]^{f_0}\ar[d]^{\varphi_0} & A_1 \ar[r]^{f_1}\ar[d]^{\varphi_1} & A_2 \ar[r]^{f_2}\ar@{-->}[d]^{\varphi_2} & \cdots\ar[r]^{f_{n-1}}&A_n\ar@{-->}[d]^{\varphi_n}\ar[r]^{f_{n}}& A_{n+1} \ar[r]^{f_{n+1}}\ar@{-->}[d]^{\varphi_{n+1}} & \Sigma^n A_0 \ar[d]^{\varepsilon_{B_{n+1}}\circ \Sigma^n \varphi_0 }\\
\Omega^nB_{n+1} \ar[r]^{\quad g_{-1}} & B_0 \ar[r]^{g_0} & B_1 \ar[r]^{g_1} & \cdots \ar[r]^{g_{n-2}\;\;}&B_{n-1}\ar[r]^{g_{n-1}}& B_{n} \ar[r]^{g_{n}}& B_{n+1}
}$$
where the upper row is in $\Theta$ and the lower row is in $\Phi$, there are
morphisms $\varphi_2,\varphi_3,\cdots,\varphi_{n+1}$ making the diagram commutative.

\item[(v)] Each solid commutative diagram in $\C$:
$$\xymatrix{
A_0 \ar[r]^{f_0}\ar[d]_{\Omega^n\varphi\circ\eta_{A_0}} & A_1 \ar[r]^{f_1}\ar@{-->}[d]^{\varphi_1} & A_2 \ar[r]^{f_2}\ar@{-->}[d]^{\varphi_2} & \cdots\ar[r]^{f_{n-1}}&A_n\ar@{-->}[d]^{\varphi_n}\ar[r]^{f_{n}}& A_{n+1} \ar[r]^{f_{n+1}}\ar[d]^{\varphi_{n+1}} & \Sigma^n A_0 \ar[d]^{\varphi}\\
\Omega^nB_{n+1} \ar[r]^{\quad g_{-1}} & B_0 \ar[r]^{g_0} & B_1 \ar[r]^{g_1} & \cdots\ar[r]^{g_{n-2}\;\;}&B_{n-1} \ar[r]^{g_{n-1}}& B_{n} \ar[r]^{g_{n}\;\;}& B_{n+1}
}$$
where the upper row is in $\Theta$ and the lower row is in $\Phi$, there are
morphisms $\varphi_1,\varphi_2,\cdots,\varphi_{n}$ making the diagram commutative.
\end{itemize}
A {\bf pre-$\bm{(n+2)}$-angulated category} is an additive category together with a pre-$(n+2)$-angulation,
and is denoted by $(\C,\Sigma^n,\Omega^n,\Theta,\Phi,\varepsilon,\eta)$.
\end{definition}

\begin{remark}
In Definition \ref{pre}, when $n=1$, it is just \cite[Definition 1.1]{BR}.
\end{remark}

We now give some examples of pre-$(n+2)$-angulated categories. From these examples,
we know that pre-$(n+2)$-angulated categories provide a common generalization of $n$-abelian,
$(n+2)$-angulated.

\begin{example}
(1) $(n+2)$-angulated categories are pre-$(n+2)$-angulated. Here $\Theta=\Phi$ and
$\Sigma^{n}=\Omega^{-n}$.

(2) We recall the notions of $n$-kernels and $n$-cokernels in an additive category from
\cite[Definiton 2.2]{Ja} and \cite[Definition 2.4]{L2}.
Let $\A$ be an additive category and $f_0\colon A_0\rightarrow A_1$  a morphism in
$\A$. An $n$-\emph{cokernel} of $f_0$ is a sequence
$$(f_1,f_2,\cdots,f_{n})\colon A_1\xrightarrow{~f_1~}A_2\xrightarrow{~f_2~}\cdots \xrightarrow
{~f_{n-1}~}A_n\xrightarrow
{~f_{n}~}A_{n+1}$$
such that the induced sequence of abelian groups
$$\xymatrix{0\xrightarrow{~~}\A(A_{n+1},B)\xrightarrow{~~} \A(A_{n},B)\xrightarrow{~~}\cdots \xrightarrow{~~}
\A(A_{1},B)\xrightarrow{~~}
\A(A_{0},B)}$$
is exact for each object $B\in\A$. We can define \emph{$n$-kernel} dually.

Any additive category with $n$-kernels and $n$-cokernels (in particular any $n$-abelian
category in the sense of Jasso \cite[Definition 3.1]{Ja}) is pre-$(n+2)$-angulated, with $\Omega^n=\Sigma^n=0$
 and $\Phi$ the class of left $n$-exact sequences
and $\Theta$ the class of right $n$-exact sequences.
\end{example}

\section{Idempotent completion of pre-$(n+2)$-angulated categories}
In this section, we first recall some basic notions on the idempotent completion of additive categories
from \cite{BM}.

Let $\C$ be an additive category. An idempotent morphism $e\colon
A \to A$ is said to be \emph{split} if there are two morphisms $p\colon A \to B$ and $q\colon B \to A$ such
that $e=qp$ and $pq=1_B$. An additive category $\C$ is said to be \emph{idempotent complet}e provided each idempotent
morphism splits.

 \begin{definition}\cite[Definition 1.2]{BM}
 Let $\C$ be an additive category. The {\em idempotent completion} of $\C$ is the category $\widetilde{\C}$ defined as follows:

  Objects of $\widetilde{\C}$ are pairs $\widetilde{A}=(A,e_{a})$, where $A$ is an object of $\C$ and $e_{a}: A\rightarrow A$ is an idempotent morphism.

  A morphism in $\widetilde{\C}$ from $(A,e_{a})$ to $(B,e_{b})$ is a morphism $p: A\rightarrow B$ in $\C$ such that $p e_{a}=e_{b}p=p$. That is to say, we have the following commutative diagram
\[
 \begin{tikzpicture}
 \draw (0,0) node{\xymatrix{
 A \ar[r]^{{e_{a}}} \ar[dr]|{p}\ar@{}[dr] \ar[d]_{p} & A \ar[d]^{p}\\
 B  \ar[r]_{e_{b}} &B.
}};
\draw (0.2,0.2) node{\tiny $\circlearrowright$};
\draw (-0.2,-0.2) node{\tiny $\circlearrowleft$};
 \end{tikzpicture}
\]
\end{definition}

Let $(\C,\Sigma^n,\Theta)$ be a right $(n+2)$-angulated category in the sense of Li \cite[Definition 2.1]{L2}.
 and $\widetilde{\C}$ be  the idempotent completion
of $\C$. Define
$\widetilde{\Sigma^n}\colon\widetilde{\C}\to\widetilde{\C}$ by $\widetilde{\Sigma}(A,e)=
(\Sigma^n A,\Sigma^ne)$. For convenience, we usually write
$\widetilde{\Sigma^n}$ as $\Sigma^n$. Define a diagram in $\widetilde{\C}$
$$A_{\bullet}\colon ~~~A_0\xrightarrow{f_0}A_1\xrightarrow{f_1}A_2\xrightarrow{f_2}\cdots\xrightarrow{f_{n-1}}A_n\xrightarrow{f_n}A_{n+1}\xrightarrow{f_{n+1}}\Sigma^n A_0$$
to be a right $(n+2)$-angle when it is a direct summand of a right $(n+2)$-angle in $\Theta$, that is,
there exists a right $(n+2)$-angle $B_{\bullet}$ in $\Theta$ and right $(n+2)$-angle morphisms
$\pi\colon A_{\bullet}\to B_{\bullet}$ and $i\colon B_{\bullet}\to A_{\bullet}$
with $\pi i=1_{B_{\bullet}}$.
Equivalently, there exists a right $(n+2)$-angle $C_{\bullet}$ in
$\widetilde{\C}$ such that $A_{\bullet}\oplus C_{\bullet}$ is
isomorphic to a right $(n+2)$-angle in $\Theta$.
We denote by $\widetilde{\Theta}$ the class of right $(n+2)$-angles in $\widetilde{\C}$.

The following theorem was
proved in \cite[Theorem 3.1]{L3}, when $\C$ is an $(n+2)$-angulated categroy. However, it can be also
extended to a right $(n+2)$-angulated categroy by Wu \cite[Theorem 3.11]{W}.

\begin{theorem}{\rm \cite[Theorem 3.1]{L3} and \cite[Theorem 3.11]{W}}\label{idem}
Let $(\C,\Sigma^n,\Theta)$ be a right $(n+2)$-angulated category.
Then Then its idempotent completion $(\widetilde{\C},\widetilde{\Sigma^n},\widetilde{\Theta})$
 admits a unique structure of right $(n+2)$-angulated category.
\end{theorem}

Note that all notions and results for right $(n+2)$-angulated categories
can be given dually for left $(n+2)$-angulated categories.
\smallskip

From now on, we assume that $(\C,\Sigma^n,\Omega^n,\Theta,\Phi,\varepsilon,\eta)$
is a pre-$(n+2)$-angulated category.

\begin{lemma}\label{key}
Each solid commutative diagram in $\widetilde{\C}$:
$$\xymatrix{
A_{\bullet}\colon &A_0 \ar[r]^{f_0}\ar[d]^{\varphi_0} & A_1 \ar[r]^{f_1}\ar[d]^{\varphi_1} & A_2 \ar[r]^{f_2}\ar@{-->}[d]^{\varphi_2} & \cdots\ar[r]^{f_{n-1}}&A_n\ar@{-->}[d]^{\varphi_n}\ar[r]^{f_{n}}& A_{n+1} \ar[r]^{f_{n+1}}\ar@{-->}[d]^{\varphi_{n+1}} & \Sigma^n A_0 \ar[d]^{\varepsilon_{B_{n+1}}\circ \Sigma^n \varphi_0 }\\
B_{\bullet}\colon&\Omega^nB_{n+1} \ar[r]^{\quad g_{-1}} & B_0 \ar[r]^{g_0} & B_1 \ar[r]^{g_1} & \cdots \ar[r]^{g_{n-2}\;\;}&B_{n-1}\ar[r]^{g_{n-1}}& B_{n} \ar[r]^{g_{n}}& B_{n+1}
}$$
where the upper row is in $\widetilde{\Theta}$ and the lower row
is in $\widetilde{\Phi}$. Then
there are
morphisms $\varphi_2,\varphi_3,\cdots,\varphi_{n+1}$ making the whole diagram commutative.
\end{lemma}

\proof By the definition of right $(n+2)$-angles in $\widetilde{\Theta}$,
there exists a right $(n+2)$-angle
$$C_{\bullet}\colon ~~~C_0\xrightarrow{h_0}C_1\xrightarrow{h_1}C_2\xrightarrow{h_2}\cdots\xrightarrow{h_{n-1}}C_n\xrightarrow{h_n}C_{n+1}\xrightarrow{h_{n+1}}\Sigma^n C_0$$
in $\Theta$ and  right $(n+2)$-angle morphisms
$$\psi_{\bullet}=(\psi_0,\psi_1,\cdots,\psi_{n+1})\colon A_{\bullet}
\to C_{\bullet}~\mbox{and}~\phi_{\bullet}=(\phi_0,\phi_1,\cdots,\phi_{n+1})\colon C_{\bullet}\to A_{\bullet}$$ such that
$\phi\psi=1_{A_{\bullet}}$.

By the definition of left $(n+2)$-angles in $\widetilde{\Phi}$,
there exists a left $(n+2)$-angle
$$D_{\bullet}\colon ~~~\Omega^n D_{n+1}\xrightarrow{u_{-1}}D_0\xrightarrow{u_0}D_1\xrightarrow{u_1}D_2\xrightarrow{u_2}
\cdots\xrightarrow{u_{n-1}}D_n\xrightarrow{u_n}D_{n+1}$$
in $\Phi$ and left $(n+2)$-angle morphisms
$$\mu_{\bullet}=(\mu_0,\mu_1,\cdots,\mu_{n+1})\colon B_{\bullet}
\to D_{\bullet}~\mbox{and}~\nu_{\bullet}=(\nu_0,\nu_1,\cdots,\nu_{n+1})\colon D_{\bullet}\to B_{\bullet}$$ such that
$\nu\mu=1_{D_{\bullet}}$.

Consider the following diagram with the first commutative square:
$$\xymatrix{
C_{\bullet}\colon &C_0 \ar[r]^{h_0}\ar[d]_{\Omega^n\mu_{n+1}\circ\varphi_0\circ\phi_0} & C_1 \ar[r]^{h_1}\ar[d]^{\mu_0\circ\varphi_1\circ\phi_1} & C_2 \ar[r]^{h_2} & \cdots\ar[r]^{h_{n-1}}&C_n\ar[r]^{h_{n}}& C_{n+1} \ar[r]^{h_{n+1}} & \Sigma^n C_0 \ar[d]_{\mu_{n+1}\circ(\varepsilon_{D_{n+1}}\circ \Sigma^n \varphi_0)\circ\Sigma^n\phi_0 }\\
D_{\bullet}\colon&\Omega^nD_{n+1} \ar[r]^{\quad u_{-1}} & D_0 \ar[r]^{u_0} & D_1 \ar[r]^{u_1} & \cdots \ar[r]^{u_{n-2}\;\;}&D_{n-1}\ar[r]^{u_{n-1}}& D_{n} \ar[r]^{u_{n}}& D_{n+1}
}$$
where the upper row is in $\Theta$ and the lower row
is in $\Phi$.  By the definition
of pre-$(n+2)$-angulation, there are morphism $\varphi'_k\colon C_k\to B_{k-1}~(k=2,3,\cdots,n+1)$
such that that makes the above diagram
commutative. We put $\varphi_k=\nu_k\varphi_k\psi_k\colon A_k\to B_{k-1}~(k=2,3,\cdots,n+1)$, as required.
This  completes the proof.  \qed

\begin{theorem}\label{main0}
Let $(\C,\Sigma^n,\Omega^n,\Theta,\Phi,\varepsilon,\eta)$ be a pre-$(n+2)$-angulated category.
Then its idempotent completion $(\widetilde{\C},\widetilde{\Sigma^n},\widetilde{\Omega^n},\widetilde{\Theta},\widetilde{\Phi},\widetilde{\varepsilon},\widetilde{\eta})$  admits a unique structure of pre-$(n+2)$-angulated category.
\end{theorem}

\proof This follows from Definition \ref{pre}, \cite[Lemma 3.1]{LS}, Theorem \ref{idem} and its dual,
Lemma \ref{key} and its dual.  \qed

\begin{remark}
In Theorem \ref{main0}, when $n=1$, it is just \cite[Theorem 3.4]{LS}.
\end{remark}

\section{Pre-$(n+2)$-angulated quotient categories}
In this section, we give more examples of pre-$(n+2)$-angulated category.
We show that how to equip quotient categories with a pre-$(n+2)$-angulated structure and give some examples to explain our main results.

\subsection{Strongly functorially finite subcategories}
Let $\C$ be  an additive category,
and $\X$ be a subcategory of $\C$.
Recall that we say a morphism $f\colon A \to B$ in $\C$ is an $\X$-\emph{monic} if
$$\Hom_{\C}(f,X)\colon \Hom_{\C}(B,X) \to \Hom_{\C}(A,X)$$
is an epimorphism for all $X\in\X$. We say that $f$ is an $\X$-\emph{epic} if
$$\Hom_{\C}(X,f)\colon \Hom_{\C}(X,A) \to \Hom_{\C}(X,B)$$
is an epimorphism for all $X\in\X$.
Similarly,
we say that $f$ is a left $\X$-approximation of $B$ if $f$ is an $\X$-monic and $A\in\X$.
We say that $f$ is a right $\X$-approximation of $A$ if $f$ is an $\X$-epic and $B\in\X$.

A subcategory $\X$ is called \emph{contravariantly finite} if any object in $\C$ admits a right
$\X$-approximation. Dually we can define  \emph{covariantly finite} subcategory.
A contravariantly finite and covariantly finite subcategory is called \emph{functorially finite}.

 We denote by $\C/\X$
the category whose objects are objects of $\C$ and whose morphisms are elements of
$\Hom_{\C}(A,B)/\X(A,B)$ for $A,B\in\C$, where $\X(A,B)$ the subgroup of $\Hom_{\C}(A,B)$ consisting of morphisms
which factor through an object in $\X$.
Such category is called the (additive) quotient category
of $\C$ by $\X$. For any morphism $f\colon A\to B$ in $\C$, we denote by $\overline{f}$ the image of $f$ under
the natural quotient functor $\C\to\C/\X$.

\begin{definition}\label{dd1}\cite[Definition 3.1]{LZ}
Let $(\C,\E,\s)$ be an $n$-exangulated category. A subcategory $\X$ of $\C$ is called
\emph{strongly contravariantly finite}, if for any object $C\in\C$, there exists a distinguished $n$-exangle
$$B\xrightarrow{}X_1\xrightarrow{}X_2\xrightarrow{}\cdots\xrightarrow{}X_{n-1}\xrightarrow{}X_{n}\xrightarrow{~g~}C\overset{}{\dashrightarrow}$$
where $g$ is a right $\X$-approximation of $C$ and $X_i\in\X$.
Dually we can define \emph{strongly covariantly finite} subcategory.

A strongly contravariantly finite and strongly  covariantly finite subcategory is called \emph{ strongly functorially finite}.
\end{definition}

\subsection{The construction of pre-$(n+2)$-angulated quotient categories}
Let $\X$ be a strongly covariantly finite subcategory of an $n$-exangulated category $\C$.
We construct a functor $\G\colon \C/\X\to\C/\X$  as follows:

For any object $C\in\C$, since $\X$ is a strongly covariantly finite subcategory of $\C$,
then there exists a distinguished $n$-exangle
\begin{equation}\label{t1}
C\xrightarrow{~a~}X_1\xrightarrow{f_1}X_2\xrightarrow{f_2}\cdots\xrightarrow{f_{n-2}}
X_{n-1}\xrightarrow{f_{n-1}}X_{n}\xrightarrow{~b~}\G C\overset{\delta}{\dashrightarrow}
\end{equation}
where $a$ is a left $\X$-approximation of $C$ and $X_1,X_2,\cdots,X_n\in\X$.
We define $\G(C)=\G C$.

For any morphism $\varphi\colon C\to C'$ in $\C$, we have the following commutative diagram
of distinguished $n$-exangles
$$\xymatrix{
C\ar[r]^{a}\ar[d]^{\varphi}& X_1 \ar[r]^{f_1}\ar@{-->}[d]^{\varphi_1} & X_2 \ar[r]^{f_2}\ar@{-->}[d]^{\varphi_2}  & \cdots \ar[r]^{f_{n-1}}& X_n \ar[r]^{b}\ar@{-->}[d]^{\varphi_n}&\G C\ar@{-->}[r]^{\delta}\ar@{-->}[d]^{\psi} & \\
C'\ar[r]^{c}&X'_1 \ar[r]^{g_1} & X'_2 \ar[r]^{g_2}  & \cdots \ar[r]^{g_{n-1}} & X'_n \ar[r]^{d}&\G C'\ar@{-->}[r]^{\delta'}&.}$$
We put $\G(\overline{\varphi})=\overline{\psi}$.
It is not difficult to see that the above assignment
induces a well-defined functor
$\G\colon \C/\X\to\C/\X$
on the quotient category $\C/\X$.

Dually
if $\X$ be a strongly contravariantly finite subcategory of an $n$-exangulated category $\C$,
then for any object $B\in\C$, there exists a distinguished $n$-exangle
\begin{equation}\label{t2}
\H B\xrightarrow{~u~}Y_1\xrightarrow{h_1}Y_2\xrightarrow{h_2}\cdots\xrightarrow{h_{n-2}}
Y_{n-1}\xrightarrow{h_{n-1}}Y_{n}\xrightarrow{~v~}B\overset{\eta}{\dashrightarrow}
\end{equation}
where $v$ is a right $\X$-approximation of $B$ and $Y_1,Y_2,\cdots,Y_n\in\X$.
We can construct a functor $\H\colon \C/\X\to\C/\X$ in a dual manner.

\begin{lemma}\label{adjoint}
Let $(\C,\E,\s)$ be an $n$-exangulated category and
$\X$ be a strongly functorially finite subcategory of $\C$. Then
there exists an adjoint pair
$$
(\G,\H) \, \colon \, \xymatrix@C=2.75pc {\C/\X \,\, \ar@<0.5ex>[r]^-{{\mathsf{}}} & \,\, \ar@<0.5ex>[l]^-{{\,\, \mathsf{}}}\C/\X}.
$$
\end{lemma}

\proof For any object $C\in\C$, using the distinguished $n$-exangles (\ref{t1}) and
(\ref{t2}),  we have the following morphism of distinguished $n$-exangles:
$$\xymatrix{
\H C\ar[r]^{a}\ar@{=}[d]& X_1 \ar[r]^{f_1}\ar[d]^{\varphi_1} & X_2 \ar[r]^{f_2}\ar[d]^{\varphi_2}  & \cdots \ar[r]^{f_{n-1}}& X_n \ar[r]^{b\;\;}\ar[d]^{\varphi_n}&\G\H C\ar@{-->}[r]^{\quad\delta}\ar[d]^{\varepsilon_C} & \\
\H C\ar[r]^{u}&Y_1 \ar[r]^{h_1} & Y_2 \ar[r]^{h_2}  & \cdots \ar[r]^{h_{n-1}} & Y_n \ar[r]^{v}&C\ar@{-->}[r]^{\eta}&}$$
and
$$\xymatrix{
C\ar[r]^{p}\ar[d]^{\eta_C}& X'_1 \ar[r]^{g_1}\ar[d]^{\phi_1} & X'_2 \ar[r]^{g_2}\ar[d]^{\phi_2}  & \cdots \ar[r]^{g_{n-1}}& X'_n \ar[r]^{q\;\;}\ar[d]^{\phi_n}&\G C\ar@{-->}[r]^{\theta}\ar@{=}[d]& \\
\H\G C\ar[r]^{s}&Y'_1 \ar[r]^{w_1} & Y'_2 \ar[r]^{w_2}  & \cdots \ar[r]^{w_{n-1}} & Y'_n \ar[r]^{t}&\G C\ar@{-->}[r]^{\omega}&}$$
where $v$ and $t$ are right $\X$-approximations, $a$ and $p$ are left $\X$-approximations.
It is straightforward to verify that the induced maps
$\overline{\varepsilon_C}\colon \G\H C\to C$ and $\overline{\eta_C}\colon C\to \H\G C$ are natural and define the counit and the
unit of an adjoint pair $(\G,\H)$ in $\C/\X$.  \qed

\medskip

\begin{lemma}\label{lem25}
Let $(\C,\E,\s)$ be an $n$-exangulated category and
$\X$ be a strongly functorially finite subcategory of $\C$. Then the following assertions hold.
\begin{itemize}
\item[\rm (1)] The quotient category $\C/\X$ is a right $(n+2)$-angulated category.

\item[\rm (2)] The quotient category $\C/\X$ is a left $(n+2)$-angulated category.

\end{itemize}
\end{lemma}

\proof (1) The proof of this lemma is just an analogue of \cite[Theorem 3.11]{ZW}.
Here we omit the proof. For the convenience of the reader, we review the specific construction as follows.
\begin{itemize}
\item  $\G\colon \C/\X\to \C/\X$ is a well-defined endofunctor in Lemma \ref{adjoint}.

\item  For any distinguished $n$-exangle
 $$A_0\xrightarrow{g_0}A_1\xrightarrow{g_1}A_2\xrightarrow{g_2}A_3
 \xrightarrow{g_3}\cdots\xrightarrow{g_{n-1}}A_n\xrightarrow{g_n}A_{n+1}\overset{\eta}{\dashrightarrow}$$
where $g_0$ is an $\X$-monic, take the following commutative diagram of distinguished $n$-exangles
$$\xymatrix{
A_0 \ar[r]^{g_0}\ar@{=}[d]& A_1 \ar[r]^{g_1}\ar@{-->}[d]^{\varphi_1} & A_2 \ar[r]^{g_2}\ar@{-->}[d]^{\varphi_2}  & \cdots \ar[r]^{g_{n-1}}& A_n \ar[r]^{g_n}\ar@{-->}[d]^{\varphi_n}&A_{n+1}\ar@{-->}[r]^{\eta}\ar@{-->}[d]^{\varphi_{n+1}} & \\
A_0 \ar[r]^{f_0}&X_1 \ar[r]^{f_1} & X_2 \ar[r]^{f_2}  & \cdots \ar[r]^{f_{n-1}} & X_n \ar[r]^{f_n}&\mathbb{G}A_0\ar@{-->}[r]^{\delta}&.}$$
Then we have a complex
$$A_0\xrightarrow{~\overline{g_0}~} A_1\xrightarrow{~\overline{g_1}~}A_2\xrightarrow
{~\overline{g_2}~}\cdots\xrightarrow{~\overline{g_{n-1}}~}A_{n}\xrightarrow{~\overline{g_{n}}~}A_{n+1}\xrightarrow{~\overline{\varphi_{n+1}}~}\mathbb{G}A_0.$$ Define right $(n+2)$-angles in $\C/\X$ as the complexes which are isomorphic to complexes obtained in this way.
We denote by $\Theta$ a class of right $(n+2)$-angles in $\C/\X$.
\end{itemize}

(2) It is dual of (1). Let $\Phi$ be a class of the left $(n+2)$-angles of $\C/\X$.
Dually we can prove $(\C/\X,\H,\Phi)$ is a left $(n+2)$-angulated category where
$\H\colon \C/\X\to \C/\X$ is a well-defined endofunctor in Lemma \ref{adjoint}.  \qed
\vspace{2mm}

The following result shows
how to equip quotient categories with a pre-$(n+2)$-angulated structure.

\begin{theorem}\label{main1}
Let $(\C,\E,\s)$ be an $n$-exangulated category and
$\X$ be a strongly functorially finite subcategory of $\C$.
Then the quotient category $\C/\X$ is a pre-$(n+2)$-angulated category.
\end{theorem}

\proof By Lemma \ref{adjoint}, we know that $(\G,\H)$ is an adjoint pair,
moveover, $\overline{\varepsilon}\colon \G\H\to {\rm Id}_{\C/\X}$
is the counit and $\overline{\eta}\colon {\rm Id}_{\C/\X}\to \H\G$
is the unit of the adjoint pair $(\G,\H)$.

By Lemma \ref{lem25}, we have that
$(\C/\X,\G,\Theta)$ is a right $(n+2)$-angulated category
and $(\C/\X,\H,\Phi)$ is a left $(n+2)$-angulated category.

So far, we have verified that $(\C/\X,\G,\H,\Theta,\Phi,\overline{\varepsilon},\overline{\eta})$ satisfies Definition \ref{pre} (i), (ii), (iii).

It remains to show that Definition \ref{pre} (iv) holds. Dually we can prove Definition \ref{pre} (v) holds.
Now we start to prove.

Given a solid commutative diagram in $\C/\X$:
$$\xymatrix{
A_0 \ar[r]^{\overline{f_0}}\ar[d]^{\overline{\varphi_0}} & A_1 \ar[r]^{\overline{f_1}}\ar[d]^{\overline{\varphi_1}} & A_2 \ar[r]^{\overline{f_2}}& \cdots\ar[r]^{\overline{f_{n-1}}}&A_n\ar[r]^{\overline{f_{n}}}& A_{n+1} \ar[r]^{\overline{f_{n+1}}}& \G A_0 \ar[d]^{\overline{\varepsilon}_{B_{n+1}}\circ \Sigma^n \overline{\varphi_0} }\\
\H B_{n+1} \ar[r]^{\quad \overline{g_{-1}}} & B_0 \ar[r]^{\overline{g_0}} & B_1 \ar[r]^{\overline{g_1}} & \cdots \ar[r]^{\overline{g_{n-2}}\;\;}&B_{n-1}\ar[r]^{\;\;\overline{g_{n-1}}}& B_{n} \ar[r]^{\overline{g_{n}}}& B_{n+1}
}$$
where the upper row is in $\Theta$ and the lower row is in $\Phi$.

For the object $\H B_{n+1}\in\C$, there exists a distinguished $n$-exangle
$$
\H B_{n+1}\xrightarrow{~x_0~}X_1\xrightarrow{x_1}X_2\xrightarrow{x_2}\cdots\xrightarrow{x_{n-2}}
X_{n-1}\xrightarrow{x_{n-1}}X_{n}\xrightarrow{~x_n~}\G\H B_{n+1}\overset{\delta}{\dashrightarrow}
$$
where $x_0$ is a left $\X$-approximation of $\H B_{n+1}$ and $X_1,X_2,\cdots,X_n\in\X$.

For the $\E$-extension $(g_{-1})_{\ast}\delta\in\E(\G\H B_{n+1},B_0)$,  there
exists a distinguished $n$-exangle
$$
B_0\xrightarrow{~c_0~}C_1\xrightarrow{c_1}C_2\xrightarrow{c_2}\cdots\xrightarrow{c_{n-2}}
C_{n-1}\xrightarrow{c_{n-1}}C_{n}\xrightarrow{~c_n~}\G\H B_{0}\overset{(g_{-1})_{\ast}\delta}{\dashrightarrow}.
$$
By (EA2$\op$), we can observe that $(g_{-1},1_{\G\H B_{0}})$ has a good lift
$\psi_{\bullet}=(g_{-1},\psi_1,\psi_2,\cdots,\psi_n, 1_{\G\H B_{0}})$. That is,
we have the morphism of distinguished $n$-exangles
$$\xymatrix{
\H B_{n+1}\ar[r]^{\quad x_0}\ar[d]^{g_{-1}}& X_1 \ar[r]^{x_1}\ar[d]^{\psi_1} & X_2 \ar[r]^{x_2}\ar[d]^{\psi_2}  & \cdots \ar[r]^{x_{n-1}}&X_n \ar[r]^{x_n\quad}\ar[d]^{\psi_n}&\G\H B_{n+1}\ar@{-->}[r]^{\qquad\delta}\ar@{=}[d] & \\
B_0 \ar[r]^{c_0}&C_1 \ar[r]^{c_1} & C_2 \ar[r]^{c_2}  & \cdots \ar[r]^{c_{n-1}} & C_n \ar[r]^{c_n\quad}&\mathbb{G}\H B_0\ar@{-->}[r]^{\quad (g_{-1})_{\ast}\delta}&}$$
and the mapping cocone
$$\H B_{n+1}\xrightarrow{\left(\begin{smallmatrix}
-x_0\\ g_{-1}
\end{smallmatrix}
\right)}X_1\oplus B_0\xrightarrow{h_1:=\left(\begin{smallmatrix}
-x_1&0\\
\psi_1&c_0
\end{smallmatrix}
\right)}X_2\oplus C_1\xrightarrow{h_2:=\left(\begin{smallmatrix}
-x_2&0\\
\psi_2&c_1
\end{smallmatrix}
\right)}\cdots$$$$\xrightarrow{h_{n-1}:=\left(\begin{smallmatrix}
-x_{n-1}&0\\
\psi_{n-1}&c_{n-2}
\end{smallmatrix}
\right)}X_n\oplus C_{n-1}\xrightarrow{\left(\begin{smallmatrix}
\psi_{n}&~c_{n-1}
\end{smallmatrix}
\right)}C_n\overset{(c_n)^{\ast}\delta}{\dashrightarrow}$$
is a distinguished $n$-exangle.

For $(1_{\H B_{n+1}},c_n)\colon(c_n)^{\ast}\delta\to\delta$, there exists a lift
$\theta_{\bullet}=(1_{\H B_{n+1}}, \theta_1,\theta_2,\cdots,\theta_n,c_n)$
by (R0). Then we
have the following commutative diagram:
$$\xymatrix@C=1.1cm{
\H B_{n+1}\ar[r]^{\left(\begin{smallmatrix}
-x_0\\ g_{-1}
\end{smallmatrix}
\right)\;\;}\ar@{=}[d]& X_1\oplus B_0 \ar[r]^{h_1}\ar[d]^{\theta_1} & X_2\oplus C_1 \ar[r]^{h_2}\ar[d]^{\theta_2}  & \cdots \ar[r]^{h_{n-1}\quad}&X_n \oplus C_{n-1}\ar[r]^{\quad\left(\begin{smallmatrix}
\theta_{n}&c_{n-1}
\end{smallmatrix}
\right)}\ar[d]^{\theta_n}&C_n\ar@{-->}[r]^{(c_n)^{\ast}\delta}\ar[d]^{c_n} & \\
\H B_{n+1} \ar[r]^{\quad x_0}&X_1 \ar[r]^{x_1} & X_2 \ar[r]^{x_2}  & \cdots \ar[r]^{x_{n-1}} & X_n \ar[r]^{x_n\quad}&\mathbb{G}\H B_{n+1}\ar@{-->}[r]^{\quad\delta}&.}$$
Since $x_0$ is a left $\X$-approximation of $\H B_{n+1}$,
then we have that $\left(\begin{smallmatrix}
-x_0\\ g_{-1}
\end{smallmatrix}
\right)\colon \H B_{n+1}\to X_1\oplus B_0$ is an $\X$-monic.
Hence we obtain that
$$\H B_{n+1}\xrightarrow{~\overline{g_{-1}}~}B_0\xrightarrow{\overline{c_0}}C_1
\xrightarrow{\overline{c_2}}C_2\cdots\xrightarrow{\overline{c_{n-1}}}C_n\xrightarrow{~\overline{c_n}~}\G\H B_{n+1}$$
is a right $(n+2)$-angle in $\C/\X$. Thus there exists a morphism of
right $(n+2)$-angles:
$$\xymatrix{
A_0 \ar[r]^{\overline{f_0}}\ar[d]^{\overline{\varphi_0}} & A_1 \ar[r]^{\overline{f_1}}\ar[d]^{\overline{\varphi_1}} & A_2 \ar[r]^{\overline{f_2}}\ar@{-->}[d]^{\overline{\psi_2}}& \cdots\ar[r]^{\overline{f_{n-1}}}&A_n\ar[r]^{\overline{f_{n}}}\ar@{-->}[d]^{\overline{\psi_n}}& A_{n+1} \ar@{-->}[d]^{
\overline{\psi_{n+1}}} \ar[r]^{\overline{f_{n+1}}}& \G A_0 \ar[d]^{\G\overline{\varphi_0}}\\
\H B_{n+1} \ar[r]^{\quad \overline{g_{-1}}} & B_0 \ar[r]^{\overline{c_0}} & C_1 \ar[r]^{\overline{c_1}} & \cdots \ar[r]^{\overline{g_{n-2}}\;\;}&C_{n-1}\ar[r]^{\;\;\overline{c_{n-1}}}& C_{n} \ar[r]^{\overline{c_{n}}\quad}&\G\H B_{n+1}.
}$$
By hypothesis, we know that
$$\xymatrix{\H B_{n+1} \ar[r]^{\quad \overline{g_{-1}}} & B_0 \ar[r]^{\overline{g_0}} & B_1 \ar[r]^{\overline{g_1}} & \cdots \ar[r]^{\overline{g_{n-2}}\;\;}&B_{n-1}\ar[r]^{\;\;\overline{g_{n-1}}}& B_{n} \ar[r]^{\overline{g_{n}}}& B_{n+1}
}$$
is a left $(n+2)$-angle in $\C/\X$, thus
we may
assume that it is induced by the following commutative diagram of distinguished $n$-exangles
$$\xymatrix{
\H B_{n+1}\ar[r]^{\quad y_0}\ar[d]^{g_{-1}}& Y_1 \ar[r]^{y_1}\ar[d]^{\phi_1} & Y_2 \ar[r]^{y_2}\ar[d]^{\phi_2}  & \cdots \ar[r]^{y_{n-1}}&Y_n \ar[r]^{y_n\;\;}\ar[d]^{\phi_n}& B_{n+1}\ar@{-->}[r]^{\quad\omega}\ar@{=}[d] & \\
B_0 \ar[r]^{g_0}&B_1 \ar[r]^{g_1} & B_2 \ar[r]^{g_2}  & \cdots \ar[r]^{g_{n-1}} & B_n \ar[r]^{g_n\;\;}&B_{n+1}\ar@{-->}[r]^{\quad (g_{-1})_{\ast}\omega}&}$$
where $y_n$ is a right $\X$-approximation of $B_{n+1}$ and $Y_1,Y_2,\cdots,Y_n\in\X$.

Since $x_0$ is a left $\X$-approximation of $\H B_{n+1}$,
there exists a morphism $\beta_1\colon X_1\to Y_1$ such that $\beta_1x_0=y_0$.
Hence we have the following commutative diagram:
\begin{equation}\label{t25}
\begin{array}{l}
\xymatrix{
\H B_{n+1}\ar[r]^{\quad x_0}\ar[d]^{g_{-1}}& X_1 \ar[r]^{x_1}\ar[d]^{\beta_1} & X_2 \ar[r]^{x_2}\ar[d]^{\beta_2}  & \cdots \ar[r]^{x_{n-1}}&X_n \ar[r]^{x_n\quad}\ar[d]^{\beta_n}&\G\H B_{n+1}\ar@{-->}[r]^{\qquad\delta}\ar[d]^{\varepsilon_{B_{n+1}}} & \\
\H B_{n+1}\ar[r]^{\quad y_0}&Y_1 \ar[r]^{y_1} & Y_2 \ar[r]^{y_2}  & \cdots \ar[r]^{y_{n-1}} & Y_n \ar[r]^{y_n\;\;}&B_{n+1}\ar@{-->}[r]^{\quad\omega}&}
\end{array}
\end{equation}
Since the diagram (\ref{t25}) is a morphism of distinguished $n$-exangles,
then we have $\delta=\varepsilon_{B_{n+1}}^\ast\omega$.

For $(1_{B_0},\varepsilon_{B_{n+1}})\colon (g_0)_{\ast}\delta\to (g_0)_{\ast}\omega$,
 there exists a lift
$\gamma_{\bullet}=(1_{B_0}, \gamma_1,\gamma_2,\cdots,\gamma_n,\varepsilon_{B_{n+1}})$
by (R0) and $(g_0)_{\ast}\delta=(g_0)_{\ast}\varepsilon_{B_{n+1}}^\ast\omega=
(\varepsilon_{B_{n+1}})_{\ast}g^{\ast}_0\omega$.
Then we have the following commutative diagram:
$$\xymatrix{
B_0\ar[r]^{c_0}\ar@{=}[d]& C_1 \ar[r]^{c_1}\ar[d]^{\gamma_1} & C_2 \ar[r]^{c_2}\ar[d]^{\gamma_2}  & \cdots \ar[r]^{c_{n-1}}&C_n \ar[r]^{c_n\quad}\ar[d]^{\gamma_n}&\G\H B_{n+1}\ar@{-->}[r]^{\qquad(g_0)_{\ast}\delta}\ar[d]^{\varepsilon_{B_{n+1}}} & \\
B_0 \ar[r]^{g_0}&B_1 \ar[r]^{g_1} & B_2 \ar[r]^{c_2}  & \cdots \ar[r]^{g_{n-1}} & B_n \ar[r]^{g_n\;\;}&B_{n+1}\ar@{-->}[r]^{\quad (g_0)_{\ast}\omega}&}$$
Now for each $2\leq k\leq n+1$, define a morphism
$\overline{\varphi_k}\colon A_k\to B_{k-1}$
by $\overline{\varphi_k}=\overline{\gamma_{k-1}}\circ\overline{\psi_{k}}$.
Then there are
morphisms $\overline{\varphi_2},\overline{\varphi_3},\cdots,\overline{\varphi_{n+1}}$ making the following diagram commutative
$$\xymatrix{
A_0 \ar[r]^{\overline{f_0}}\ar[d]^{\overline{\varphi_0}} & A_1 \ar[r]^{\overline{f_1}}\ar[d]^{\overline{\varphi_1}} & A_2\ar@{-->}[d]^{\overline{\varphi_2}} \ar[r]^{\overline{f_2}}& \cdots\ar[r]^{\overline{f_{n-1}}}&A_n\ar@{-->}[d]^{\overline{\varphi_n}} \ar[r]^{\overline{f_{n}}}& A_{n+1}\ar@{-->}[d]^{\overline{\varphi_{n+1}}}  \ar[r]^{\overline{f_{n+1}}}& \G A_0 \ar[d]^{\overline{\varepsilon}_{B_{n+1}}\circ \Sigma^n \overline{\varphi_0} }\\
\H B_{n+1} \ar[r]^{\quad \overline{g_{-1}}} & B_0 \ar[r]^{\overline{g_0}} & B_1 \ar[r]^{\overline{g_1}} & \cdots \ar[r]^{\overline{g_{n-2}}\;\;}&B_{n-1}\ar[r]^{\;\;\overline{g_{n-1}}}& B_{n} \ar[r]^{\overline{g_{n}}}& B_{n+1}.
}$$
This completes the proof.   \qed

\medskip

By applying Theorem \ref{main1} to $(n+2)$-angulated categories, and using the fact that any $(n+2)$-angulated category can be viewed as an $n$-exangulated category, we get the following result.

\begin{corollary}\label{cor1}
Let $\C$ be an $(n+2)$-angulated category and
$\X$ be a strongly functorially finite subcategory of $\C$.
Then the quotient category $\C/\X$ is a pre-$(n+2)$-angulated category.
\end{corollary}

\begin{remark}
In Corollary \ref{cor1}, when $n=1$, it is just the Theorem 2.2 in \cite{J}.
\end{remark}

By applying Theorem \ref{main1} to $n$-exact categories, and using the fact that any $n$-exact category can be viewed as an $n$-exangulated category, we get the following result.

\begin{corollary}\label{cor2}
Let $\C$ be an $n$-exact category and
$\X$ be a strongly functorially finite subcategory of $\C$.
Then the quotient category $\C/\X$ is a pre-$(n+2)$-angulated category.
\end{corollary}

Now we give some concrete examples.

\begin{example}
This example comes from \cite{L1}.
Let $\T=D^b(kQ)/\tau^{-1}[1]$ be the cluster category of type $A_3$,
where $Q$ is the quiver $1\xrightarrow{~\alpha~}2\xrightarrow{~\beta~}3$,
$D^b(kQ)$ is the bounded derived category of finite generated modules over $kQ$, $\tau$ is the Auslander-Reiten translation and $[1]$ is the shift functor
of $D^b(kQ)$. Then $\T$ is a $2$-Calabi-Yau triangulated category. Its shift functor is also denoted by $[1]$.

We describe the
Auslander-Reiten quiver of $\T$ in the following:
$$
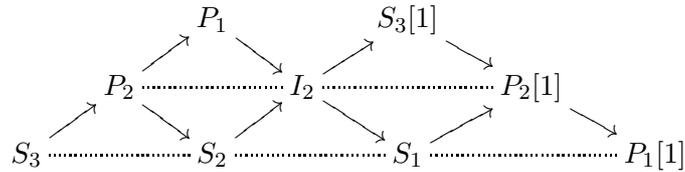
\xymatrix@C=0.6cm@R0.3cm{
&&P_1\ar[dr]
&&S_3[1]\ar[dr]
&&\\
&P_2 \ar@{.}[rr] \ar[dr] \ar[ur]
&&I_2 \ar@{.}[rr] \ar[dr] \ar[ur]
&&P_2[1]\ar[dr]\\
S_3\ar[ur]\ar@{.}[rr]&&S_2\ar[ur]\ar@{.}[rr]
&&S_1\ar[ur]\ar@{.}[rr]&&P_1[1]
}
$$
It is straightforward to verify that $\C:=\add(S_3\oplus P_1\oplus S_1)$ is a $2$-cluster tilting subcategory of $\T$ and satisfies $\C[2]=\C$.   By \cite[Theorem 1]{GKO}, we know that $\C$ is a $4$-angulated category with an automorphism functor $[2]$.
Let $\X=\add(S_3\oplus S_1)$.
Then the $4$-angle
$$P_1\xrightarrow{~~}S_1\xrightarrow{~~}S_3\xrightarrow{~~}P_1\xrightarrow{~~}P_1[2]$$
shows that $\X$ is strongly functorially finite subcategory of $\C$.
By Corollary \ref{cor1}, we have that the quotient category $\C/\X$ is a pre-$4$-angulated category.
\end{example}

\begin{example}
Let $\Lambda$ be the algebra given by the following (infinity) quiver with relations $x^4=0$:
$$\xymatrix{1 &2\ar[l]_{x} &3\ar[l]_{x} &4\ar[l]_{x} &\cdots \ar[l]_{x} &n\ar[l]_{x} &\cdots \ar[l]_{x}}$$
The Auslander-Reiten quiver of ${\rm mod}\Lambda$ is the following. {\small
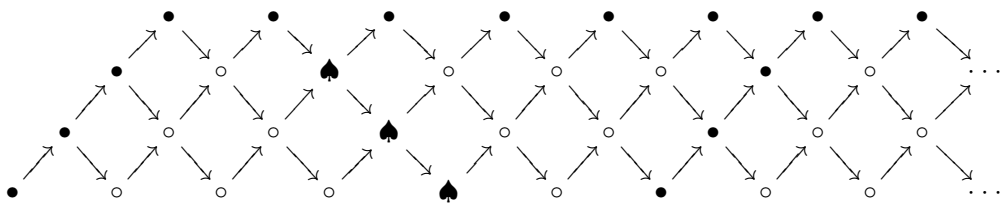
$$\xymatrix@C=0.3cm@R0.3cm{
 &&&\bullet\ar[dr] &&\bullet \ar[dr] &&\bullet \ar[dr] &&\bullet \ar[dr] &&\bullet \ar[dr] &&\bullet \ar[dr] &&\bullet \ar[dr] &&\bullet \ar[dr]\\
&&\bullet \ar[ur] \ar[dr] &&\circ \ar[ur] \ar[dr] &&\spadesuit \ar[ur]  \ar[dr] &&\circ \ar[ur] \ar[dr]  &&\circ \ar[ur] \ar[dr] &&\circ \ar[ur]  \ar[dr] &&\bullet\ar[ur] \ar[dr] &&\circ\ar[ur] \ar[dr] &&\cdot\cdot\cdot\\
 &\bullet \ar[ur] \ar[dr] &&\circ \ar[ur] \ar[dr] &&\circ \ar[ur] \ar[dr]  &&\spadesuit \ar[ur] \ar[dr] &&\circ \ar[ur] \ar[dr] &&\circ \ar[ur] \ar[dr] &&\bullet \ar[ur] \ar[dr] &&\circ \ar[ur] \ar[dr] &&\circ \ar[ur] \ar[dr]\\
\bullet \ar[ur] &&\circ \ar[ur] &&\circ \ar[ur] &&\circ \ar[ur] &&\spadesuit \ar[ur] &&\circ \ar[ur] &&\bullet \ar[ur] &&\circ \ar[ur] &&\circ \ar[ur] &&\cdot\cdot\cdot
}
$$}
where the object denoted by $\spadesuit$ and $\bullet$ appear periodically. Let $\C$ be the additive closure of all the indecomposable objects denoted by $\spadesuit$ and $\bullet$. Then $\C$ is a cluster tilting subcategory of ${\rm mod}\Lambda$, hence it is $2$-abelian (which is $2$-exact). Let $\X$ be the additive closure of all the indecomposable objects denoted by $\bullet$. Then $\X$ is a strongly functorially finite subcategory of $\C$. By Corollary \ref{cor2}, we know that $\C/\X$ is a pre-$4$-angulated category.
\end{example}

\section{Application to $(n+2)$-angulated categories}
Note that any $(n+2)$-angulated category can be viewed as an $n$-exangulated category.
We will give an application of Theorem \ref{main1}.
Unless otherwise specified $k$ will be a field and
all categories in this section will be $k$-linear Hom-finite Krull-Schmidt.

We assume that $(\C,\Sigma^n,\Theta)$ is an $(n+2)$-angulated category
and $\X$ is a strongly functorially finite subcategory of $\C$.

Recall that a morphism $f\colon A\to B$ is called \emph{$\X$-ghost} if
$\C(\X, f)=0$. Dually $f$ is called \emph{$\X$-coghost} if $\C(f,\X)=0$.
We denote by $\gh_{\X}(A, B)$ (resp. $\cogh_{\X}(A, B))$ is a subgroup $\C(A, B)$ consisting of all $\X$-ghost
 (resp. $\X$-coghost) morphisms. For more details, see \cite[Section 3.1]{B}.

Since $\X$ is strongly functorially finite, for any object $C\in\C$,
there are distinguished $n$-exangles
$$
\H C\xrightarrow{~u~}Y_1\xrightarrow{h_1}Y_2\xrightarrow{h_2}\cdots\xrightarrow{h_{n-2}}
Y_{n-1}\xrightarrow{h_{n-1}}Y_{n}\xrightarrow{~v~}C\xrightarrow{~w~}\Sigma^n\H C
$$
and
$$
C\xrightarrow{~a~}X_1\xrightarrow{f_1}X_2\xrightarrow{f_2}\cdots\xrightarrow{f_{n-2}}
X_{n-1}\xrightarrow{f_{n-1}}X_{n}\xrightarrow{~b~}\G C\xrightarrow{~c~}\Sigma^nC
$$
where $a$ is a left $\X$-approximation of $C$ and $X_1,X_2,\cdots,X_n\in\X$, and
 $v$ is a right $\X$-approximation of $C$ and $Y_1,Y_2,\cdots,Y_n\in\X$.

It is obvious that $w\colon C\to \Sigma^n\H C$ is $\X$-ghost (since $v$ is a right $\X$-approximation and $wv=0$)
and any $\X$-ghost morphism
$C\to B$ factors through $w$.
\medskip

The following result is a higher analog of \cite[Proposition 3.1]{B}.

\begin{proposition}\label{prop}
Let $(\C,\Sigma^n,\Theta)$ be an $(n+2)$-angulated category
and $\X$ be a strongly functorially finite subcategory of $\C$.
For any object $C\in\C$, the following statements hold.
\begin{itemize}
\item[\rm (i)] $\gh_{\X}(C,\Sigma^n\X)=0$ if and only if
$\G\H C\longrightarrow C$  is invertible in $\C/\X$.

\item[\rm (ii)] $\cogh_{\X}(\Sigma^{-n}\X,C)=0$ if and only if
$C\longrightarrow\H\G C$  is invertible in $\C/\X$.

\item[\rm (iii)] The adjoint pair $(\G,\H)$ induces an equivalence:
$$\{C\in\C\mid \gh_{\X}(C,\Sigma^n\X)=0\}/\X\xrightarrow{~\simeq~}
\cogh_{\X}(\Sigma^{-n}\X,C)=0/\X.$$
In particular the pre-$(n+2)$-angulated category $\C/\X$ is $(n+2)$-angulated if and only if for any object
$C\in\C$:
$$\gh_{\X}(C,\Sigma^n\X)=0=\cogh_{\X}(\Sigma^{-n}\X,C).$$
\end{itemize}
\end{proposition}

\proof  This is an adaptation of the proof of  \cite[Proposition 3.1]{B}.

 (i) ``$\Longrightarrow$"
Suppose $\gh_{\X}(C,\Sigma^n\X)=0$.
We claim that $u$ is a left $\X$-approximation of $\H C$.
In fact, for any morphism $d\colon \H C\to X$ with $X\in\X$,
since $v$ is an $\X$-ghost,
then the composition $C\xrightarrow{~w~}\Sigma^n\H C\xrightarrow{~\Sigma^nd~}\Sigma^nX$
is zero. Thus there exists a morphism $q\colon \Sigma^nY_1\to\Sigma^nX$ such that
$\Sigma^nd=q\circ (-1)^n\Sigma^nu$.
$$\xymatrix@C=1.5cm{C\ar[rd]_0\ar[r]^{w\quad}&\Sigma^n\H C\ar[d]^{\Sigma^nd}\ar[r]^{(-1)^n\Sigma^nu}&\Sigma^n Y_1\ar@{-->}[dl]^{q}\\
&\Sigma^nX&}$$
It follows that $d\colon \H C\to X$ factors through $Y_1\in\X$.
This shows that $u$ is a left $\X$-approximation of $\H C$.
By construction, see the proof of Lemma \ref{adjoint}, this clearly implies that the counit
$\overline{\varepsilon}_C\colon \G\H C\to C$ is invertible in $\C/\X$.

``$\Longleftarrow$"
We assume that the counit $\overline{\varepsilon}_C\colon\G\H C\to C$ is invertible in $\C/\X$, we prove that
the morphism $u$ is a left $\X$-approximation of $\H C$.

Since $\overline{\varepsilon}_A\colon\G\H C\to C$ is invertible in $\C/\X$, then there exists $\overline{\beta}_{C}\colon C\to \G\H C$
such that $\overline{\varepsilon}_C\circ\overline{\beta}_{C}=\overline{1}$.
Thus the morphism $1-\varepsilon_C\beta_C$ factors through some object in $\X$ and then
it factors through the right $\X$-approximation $v\colon Y_n\to C$ of $C$.
That is, there exists a morphism $p\colon \G\H C\to Y_n$ such that
$1-\varepsilon_C\beta_C=vp$. It follows that
$w(1-\varepsilon_C\beta_C)=(wv)p=0$ and then $w=w\varepsilon_C\beta_C$.

Note that we have the following commutative diagram:
$$\xymatrix{
\H C\ar[r]^{a}\ar@{=}[d]& X_1 \ar[r]^{f_1}\ar[d]^{\varphi_1} & X_2 \ar[r]^{f_2}\ar[d]^{\varphi_2}  & \cdots \ar[r]^{f_{n-1}}& X_n \ar[r]^{b\;\;}\ar[d]^{\varphi_n}&\G\H C\ar[r]^{e}\ar[d]^{\varepsilon_C} & \Sigma^n\H C\ar@{=}[d]\\
\H C\ar[r]^{u}&Y_1 \ar[r]^{h_1} & Y_2 \ar[r]^{h_2}  & \cdots \ar[r]^{h_{n-1}} & Y_n \ar[r]^{v}&C\ar[r]^{w\quad}&\Sigma^n\H C}$$
where $a$ is a left $\X$-approximation of $\H C$ and $X_1,X_2,\cdots,X_n\in\X$.
By the last square is commutative, we have $w=e\beta_C$ and then $w=w\varepsilon_C\beta_C=e\beta_C$.
Hence we have a morphism of $(n+2)$-angles
$$\xymatrix{
\H C\ar[r]^{u}\ar@{=}[d]& Y_1 \ar[r]^{h_1}\ar@{-->}[d]^{\phi} & Y_2 \ar[r]^{h_2}\ar@{-->}[d]  & \cdots \ar[r]^{h_{n-1}}& Y_n \ar[r]^{u\;\;}\ar@{-->}[d]&\G\H C\ar[r]^{w}\ar[d]^{\beta_C} & \Sigma^n\H C\ar@{=}[d]\\
\H C\ar[r]^{a}&X_1 \ar[r]^{f_1} & X_2 \ar[r]^{f_2}  & \cdots \ar[r]^{f_{n-1}} & X_n \ar[r]^{b\quad}&\G\H C\ar[r]^{e}&\Sigma^n\H C.}$$
By the first square is commutative, we have $a=\phi u$.

For any morphism $\gamma\colon \H C\to X$ with $X\in\X$,
since $a$ is left $\X$-approximation, then there exists a morphism $\mu\colon X_1\to X$
such that $\gamma=\mu a$.
Therefore $\gamma=\mu a=(\mu\phi)u$. This shows that $u$ is a left $\X$-approximation of $\H C$.

Now we prove $\gh_{\X}(C,\Sigma^n\X)=0$. Let $\psi\colon C\to \Sigma^n X$ be an $\X$-ghost morphism
where $X\in\X$. It is obvious that $\psi$ factors through $w$, that is,
there exists a morphism $\alpha\colon \Sigma^n\H C\to \Sigma^n X$ such that
$\psi=\alpha w$.
Since  $u$ is a left $\X$-approximation of $\H C$, there exists a morphism
$\varrho\colon Y_1\to X$ such that $\Sigma^{-n}\alpha=\varrho u$.
It follows that $$\Sigma^{-n}\psi=\Sigma^{-n}\alpha\circ\Sigma^{-n}w=\varrho\circ (u\circ\Sigma^{-n}w)=0$$
and then $\psi=0$.
This shows that $\gh_{\X}(C,\Sigma^n\X)=0$.

(ii) It is dual of (i).

(iii)  It follow directly from (i) and (ii).   \qed
\medskip

\begin{remark}
In Proposition \ref{prop}, when $n=1$, it is just \cite[Proposition 3.1]{B}.
\end{remark}

\subsection{Auslander-Reiten $(n+2)$-angles and Serre functors}
Let  $(\C,\Sigma^n,\Theta)$ be an $(n+2)$-angulated category. We denote by ${\rm rad}_{\C}$ the Jacobson radical of $\C$. Namely, it is given by the formula
$${\rm rad}_{\C}(X,Y)=\{g\colon X\to Y~|~{\rm Id_{\emph{X}}-}hg~\mbox{is invertible for any}~ h\colon Y\to X\},$$
for all objects $X$ and $Y$ in $\C$.

Next, we recall some terminology for the Auslander--Reiten
theory. Let $f\colon X \to Y$ be a morphism. One says that $f$ is \emph{right almost split} if $f$ is not a split epimorphism and every non-split epimorphism  $g\colon M \to  Y$ factors through $f$.
In a dual manner, one defines $f$ to be \emph{left almost split}.

\begin{definition}\cite[Definition 3.8]{IY} and \cite[Definition 5.1]{F}
Let  $(\C,\Sigma^n,\Theta)$ be an $(n+2)$-angulated category.
An $(n+2)$-angle
$$\xymatrix {A_0 \xrightarrow{~\alpha_0~}A_1 \xrightarrow{~\alpha_1~} A_2 \xrightarrow{~\alpha_2~} \cdots
  \xrightarrow{~\alpha_{n - 1}~} A_n \xrightarrow{~\alpha_{n}~} A_{n+1}\xrightarrow{~\alpha_{n+1}~} \Sigma^n A_0}$$
in $\C$ is called an \emph{Auslander-Reiten $(n+2)$-angle }if
$\alpha_0$ is left almost split, $\alpha_n$ is right almost split and
when $n\geq 2$, also $\alpha_1,\alpha_2,\cdots,\alpha_{n-1}$ are in $\rad_{\C}$.
\end{definition}

Let $\C$ be a $k$-linear Hom-finite additive category where $k$ is a field. Recall that $\C$ has a \emph{Serre functor} $\mathbb{S}$, that is, an auto-equivalence for which there are a natural equivalence
$$\C(X,Y)\simeq D\C(Y, \mathbb{S}X)$$
for any $X,Y\in\C$, where $D=\Hom_k (-, k)$ is the $k$-linear duality functor.
\medskip

The following theorem is due to \cite[Theorem 4.5]{Z}, which builds a link between Auslander-Reiten $(n+2)$-angles and Serre functors.

\begin{theorem}\label{serre}{\rm \cite[Theorem 4.5]{Z}}
Let  $(\C,\Sigma^n,\Theta)$ be an $(n+2)$-angulated category.
Then $\C$ has Auslander-Reiten $(n+2)$-angles if and only if
$\C$ has a Serre functor $\mathbb{S}$.

If either of these properties holds, then the action of the Serre functor on objects coincides with $\tau_n\Sigma^n$, namely
$\tau_n=\mathbb{S}\Sigma^{-n}$. In this case, $\tau_n$ is called  the $n$-Auslander-Reiten translation.
\end{theorem}

\subsection{$(n+2)$-angulated quotient categories}

The following result shows
how to equip quotient categories with an $(n+2)$-angulated structure, which recovers and gives a short proof to a result of Zhou \cite[Theorem 5.7]{Z}.

\begin{theorem}\label{main6}
Let $(\C,\Sigma^n,\Theta)$ be an $(n+2)$-angulated category with
 a Serre functor $\mathbb{S}$, and $\X$ be a strongly functorially finite subcategory of $\C$.
 Then the following statements hold.
 \begin{itemize}
\item[\rm (i)] $\gh_{\X}(-,\Sigma^n\X)=0$ if and only if $\Sigma^n\X\subseteq\mathbb{S}\X$.

\item[\rm (ii)] $\cogh_{\X}(\Sigma^{-n}\X,-)=0$ if and only if
$\mathbb{S}\X\subseteq\Sigma^n\X$.

\item[\rm (iii)] The pre-$(n+2)$-angulated category $\C/\X$ is $(n+2)$-angulated if and only $\mathbb{S}\X=\Sigma^n\X$.
\end{itemize}
 \end{theorem}

\proof This is an adaptation of the proof of  \cite[Corollary 3.6]{B}.

(i) ``$\Longleftarrow$"
Assume that $\beta\in\gh_{\X}(C,\Sigma^n X)$ is any $\X$-ghost where $X\in\X$.
Since $\Sigma^n\X\subseteq\mathbb{S}\X$, then there exists an object $X'\in\X$ such that
$\mathbb{S}X'=\Sigma^nX$.
Then we have that $\C(\X,\beta)\colon\C(\X,C)\to\C(\X,\Sigma^nX)$ is zero since
$\beta$ is $\X$-ghost.
Using a Serre functor $\mathbb{S}$, we obtain that the morphism
$\C(X',\X)\to\C(\mathbb{S}^{-1}C,\X)$ is also zero.
Thus The image of the identity $1_{X'}$ under this morphism is the morphism
$\mathbb{S}^{-1}\beta\colon \mathbb{S}^{-1}C\to X'$.
Therefore $\mathbb{S}^{-1}\beta=0$ implies $\beta=0$. This is what we expect.

``$\Longleftarrow$"  Assume that $0\neq X$ is any indecomposable object in $\X$.
By Theorem \ref{serre}, there exists an Auslander-Reiten $(n+2)$-angles
$$X\xrightarrow{~\alpha~}A_1\to A_2\to\cdots\to A_{n}\xrightarrow{~\beta~}\mathbb{S}^{-1}\Sigma^nX\xrightarrow{~\gamma~}\Sigma^nX.$$
This shows that $\beta$ is right almost split and $\gamma\neq 0$.

If $\mathbb{S}^{-1}\Sigma^nX\notin\X$, then  any morphism $g\colon X'\to \mathbb{S}^{-1}\Sigma^nX$ where $X'\in\X$
is not split epimorphism.
Thus it factors through $\beta$, that is, there exists a morphism $b\colon X'\to A_n$ such that
$g=\beta b$. It follows that $\gamma g=(\gamma\beta)b=0.$
This shows that the morphism $\gamma\colon \mathbb{S}^{-1}\Sigma^nX\to\Sigma^nX$ is an $\X$-ghost and therefore it is zero. This is impossible since $\gamma\neq 0$.
Therefore $\mathbb{S}^{-1}\Sigma^nX\in\X$ implies $\Sigma^nX\in\mathbb{S}\X$ for any indecomposable object $X\in\X$.
Since $\C$ is a Krull-Schmidt category we obtain that $\Sigma^n\X\subseteq\mathbb{S}\X$.

(ii) It is dual of (i).

(iii) It follows from (i), (ii) and Proposition \ref{prop}.  \qed

\begin{remark}
In Theorem \ref{main6}, when $n=1$, it is just \cite[Corollary 3.6]{B} and \cite[Theorem 3.3]{J}.
\end{remark}

\textbf{Jian He}\\
College of Science, Hunan University of Technology and Business, Changsha 410205, Hunan P. R. China\\
E-mail: \textsf{jinghe1003@163.com}\\[0.3cm]
\textbf{Panyue Zhou}\\
College of Mathematics, Hunan Institute of Science and Technology, Yueyang 414006, Hunan, P. R. China.\\
E-mail: \textsf{panyuezhou@163.com}
\\[0.3cm]
\textbf{Xingjia Zhou}\\
College of Mathematics, Hunan Institute of Science and Technology, Yueyang 414006, Hunan, P. R. China.\\
E-mail: \textsf{xingjiazhou@126.com}


\begin{thebibliography}{99}
\bibitem[ABM]{ABM} I. Assem, A. Beligiannis, N. Marmaridis.
Right triangulated categories with right semi-equivalences. Algebras and modules, II (Geiranger, 1996), 17--37, CMS Conf. Proc., 24, Amer. Math. Soc., Providence, RI, 1998.

\bibitem[B]{B} A. Beligiannis.
Rigid objects, triangulated subfactors and abelian localizations. Math. Z. 274 (3-4) (2013),  841--883.

\bibitem[BM]{BM} P. Balmer, M. Schlichting. Idempotent completion of triangulated categories. J. Algebra 236 (2001), 819--834.

\bibitem[BR]{BR} A. Beligiannis, I. Reiten. Homological and homotopical aspects of torsion theories. Mem. Amer. Math. Soc. 188 (2007), no. 883, viii+207 pp.


\bibitem[F]{F}
F. Fedele. Auslander-Reiten $(d+2)$-angles in subcategories and a $(d+2)$-angulated generalisation of a theorem by Br\"{u}ning. J. Pure Appl. Algebra, 223(8) (2019),  3554--3580.


\bibitem[GKO]{GKO} C. Geiss, B. Keller, S. Oppermann. $n$-angulated categories. J. Reine Angew. Math. 675 (2013), 101--120.


\bibitem[HHZ]{HHZ} J. He, J. He and P. Zhou. Idempotent completion of certain $n$-exangulated categories.
arXiv: 2207.01232, 2022.

\bibitem[HLN1]{HLN1} M. Herschend, Y. Liu, H. Nakaoka. $n$-exangulated categories (I): Definitions and fundamental properties. J. Algebra 570 (2021), 531--586.

\bibitem[HLN2]{HLN2} M. Herschend, Y. Liu, H. Nakaoka.  $n$-exangulated categories (II): Constructions from $n$-cluster tilting subcategories. J. Algebra 594 (2022), 636--684.

\bibitem[HZZ]{HZZ} J. Hu, D. Zhang, P. Zhou. Proper classes and Gorensteinness in extriangulated categories. J. Algebra 551 (2020), 23--60.

\bibitem[HZZ1]{HZZ1} J. Hu, D. Zhang, P. Zhou. Two new classes of $n$-exangulated categories. J. Algebra 568 (2021), 1--21.


\bibitem[IY]{IY} O. Iyama, Y. Yoshino. Mutations in triangulated categories and rigid Cohen-Macaulay modules. Invent. Math. 172(1) (2008), 117--168.

\bibitem[Ja]{Ja} G. Jasso.  $n$-abelian and $n$-exact categories. Math. Z. 283(3-4) (2016), 703--759.

\bibitem[J]{J} P. J{\o}rgensen. Quotients of cluster categories. Proc. Roy. Soc. Edinburgh Sect. A,  140(1) (2010), 65--81.


\bibitem[KMS]{KMS} C. Klapproth, D. Msapato, A. Shah.
 Idempotent completions of $n$-exangulated categories. arXiv: 2207.04023, 2022.

\bibitem[L1]{L1} Z. Lin. $n$-angulated quotient categories induced by mutation pairs.
Czechoslovak Math. J. 65(4) (2015),  953--968.

\bibitem[L2]{L2}  Z. Lin. Right $n$-angulated categories arising from covariantly finite
subcategories. Comm. Algebra, 45(2) (2017), 828--840.

\bibitem[L3]{L3}  Z. Lin. Idempotent completion of $n$-angulated categories. Appl. Categ. Structures 29(6) (2021), 1063--1071.

\bibitem[LS]{LS} J. Liu, L. Sun. Idempotent completion of pretriangulated categories.  Czechoslovak Math. J. 64(2) (2014), 477--494.



\bibitem[LZ]{LZ} Y. Liu, P. Zhou.  Frobenius $n$-exangulated categories.
J. Algebra 559 (2020), 161--183.


\bibitem[NP]{NP} H. Nakaoka, Y. Palu.  Extriangulated categories, Hovey twin cotorsion pairs and model structures. Cah. Topol. G\'{e}om. Diff\'{e}r. Cat\'{e}g. 60(2) (2019),  117--193.

\bibitem[NP1]{NP1} H. Nakaoka, Y. Palu. External triangulation of the homotopy category of exact quasi-category.
arXiv: 2004.02479, 2020.


\bibitem[V]{V} J. Verdier.
Des cat\'{e}gories d\'{e}riv\'{e}es des cat\'{e}gories ab\'{e}liennes.
With a preface by Luc Illusie. Edited and with a note by Georges Maltsiniotis.
Ast\'{e}risque No. 239 (1996).


\bibitem[W]{W} S. Wu. Idempotent completion of left
$n$-angulated categories and $n$-exact categories.  Master's thesis, 2020.


\bibitem[Z]{Z} P. Zhou. Higher--dimensional Auslander--Reiten theory on $(d+2)$-angulated categories. To appear in Glasgow Mathematical Journal, First View, pp. 1--21. DOI: https://doi.org/10.1017/S0017089521000343



\bibitem[ZW]{ZW}
Q. Zheng, J. Wei. $(n+2)$-angulated quotient categories.
 Algebra Colloq. 26(4) (2019),  689--720.

\bibitem[ZZ]{ZZ} P. Zhou, B. Zhu. Triangulated quotient categories revisited. J. Algebra 502 (2018), 196--232.

\bibitem[ZhZ]{ZhZ} B. Zhu, X. Zhuang. Tilting subcategories in extriangulated categories. Front. Math. China 15(1) (2020),  225--253.

\end{thebibliography}
\end{document}